\documentclass[11pt,reqno]{amsart}
\usepackage{amsmath,latexsym,amssymb,verbatim,amsbsy,amsthm,mathrsfs,tikz,times}
\usepackage[margin=1in]{geometry}
\usepackage{color}
\usepackage[colorlinks=true, pdfstartview=FitV, linkcolor=blue, citecolor=blue,
urlcolor=blue]{hyperref}

\theoremstyle{plain}
\newtheorem{theorem}{Theorem}[section]

\newtheorem{proposition}[theorem]{Proposition}

\theoremstyle{definition}

\def\tilde{\widetilde}
\numberwithin{equation}{section}

\renewcommand\hat{\widehat}

\def\RR{{\mathbb R}}
\def\TT{{\mathbb T}}

\renewcommand{\phi}{\varphi}

\title{\bf Improved a priori bounds for thermal fluid equations}
 \author{Andrei Tarfulea}
\address{Department of Mathematics, University of Chicago}
\email{atarfulea@math.uchicago.edu}
\date{}

\begin{document}

\maketitle
\begin{abstract}
We consider two hydrodynamic model problems (one incompressible and one compressible)
with three dimensional fluid flow on
the torus and temperature-dependent viscosity and conductivity.
The ambient heat for the
fluid is transported by the flow and fed by the local energy dissipation, modeling the
transfer of
kinetic energy into thermal energy through fluid friction. Both the viscosity and
conductivity grow with the local temperature. We
prove a strong a priori bound on the enstrophy of the velocity weighed against the temperature
for initial data of arbitrary size, requiring only that the conductivity be proportionately
larger than the viscosity (and, in the incompressible case, a bound on the temperature as
a Muckenhoupt weight).
\end{abstract}

\section{Introduction and Main Results}
\noindent{}The existence of solutions for models of fluid equations taking temperature
and density into account is a delicate problem. When the dissipation has a
local dependence on the temperature, many standard techniques for analyzing fluid equations
(e.g., a vorticity formulation or methods using paraproduct decomposition) do not readily apply.
For a general treatment of such models, see \cite{Feireisl} and \cite{Feireisl-again}.
Even establishing the physically expected thermodynamic properties,
such as (for example) a uniform lower bound on the temperature, can take great
effort; see \cite{B-V} and \cite{M-V}.
To the best of our knowledge, most investigations of such fluid equations
show existence of weak (and admissible) solutions or essentially show that the
temperature dependence does not create exotic pathologies compared with the usual
Navier-Stokes equations \cite{BFM}. Many models (including \cite{B-V} and \cite{M-V})
assume from the beginning that the viscosity and conductivity are, as functions of the
temperature, bounded above and below (see (1.8) of \cite{F-M}). They require this to
establish a minimal coercivity for the dissipation.\\
\\
In this article, however, we examine two models where
the viscosity and conductivity are directly proportional to the square root
of the temperature. This is motivated largely by formula (15.3) of \cite{Lautrup},
the empirical rate at which viscosity increases with temperature in
gases. The enhanced Brownian motion causes molecules from one portion of a gas to spread
faster throughout the medium, so that the averaged velocities (macroscopic flows) are seen
to equalize faster in hotter areas. We first consider the following coupled system of nonlinear
PDE arising from the Navier-Stokes Equations with temperature-dependent viscosity and
conductivity
\begin{align}
& \partial_t u + u \cdot \nabla u + \nabla p - \nu \text{div}(\theta \nabla u) = 0
 \ , \ \ \ \ \ \ \ \ \ \nabla \cdot u = 0 \label{thNSE-1} \\
& \partial_t (\theta^2) + u \cdot \nabla (\theta^2) -
\kappa \text{div}(\theta \nabla (\theta^2)) = \nu \theta | \nabla u |^2   \label{Temp-1}
\end{align}
on the three-dimensional torus $\TT^3$ with given (sufficiently regular)
initial data $u_0$ and $\theta_0$.
We require that $\theta_0(x) \geq 1$ for all $x \in \TT^3$ and that $u_0$ be mean-zero;
\eqref{thNSE-1} then forces any smooth solution $u$ to remain mean-zero for all time.
The constant parameters $\nu$ and $\kappa$ respectively indicate the strength of the
kinematic viscosity and heat conductivity of the fluid. We define
$$K = \kappa / \nu + 1 $$
to be (a shift of) the inverse Prandtl number; this parameter will appear many times
in the calculations that follow. We remark here that both $\nu$ and $\kappa$ should have
units of $\text{K}^{-1/2} \text{L}^{1/2} \text{T}^{-1}$ (here K is absolute temperature,
L is units of length, and
T is units of time). Modulo the common factor of $\theta$ (which has units of
$\text{K}^{1/2} \text{L}^{3/2}$ that must be canceled), $\nu$ and $\kappa$ behave
like typical ``renormalized dissipations'' (with units of $\text{L}^2 \text{T}^{-1}$).
Their ratio is therefore
unitless (as is the Prandtl number in general), and so is $K$. This is important to note
since the main results will include conditions on the size of $K$ and also use $K$
as an exponent.\\
\\
Here we use $\theta^2$ to denote the local temperature of
the fluid. However, it will be convenient from now on to work with $\theta$, the positive square
root of the temperature (see Proposition \ref{prop-1}); \eqref{Temp-1} then reduces to
\begin{equation}
\partial_t \theta + u \cdot \nabla \theta - \kappa \theta \Delta \theta -
2 \kappa | \nabla \theta |^2 = \frac{\nu}{2} | \nabla u |^2
\label{Thet-1}
\end{equation}
\noindent{}Note that the pressure term $p$ appearing in \eqref{thNSE-1} has a
different explicit representation than for the standard Navier-Stokes Equations.
Specifically,
\begin{equation}
p = (-\Delta)^{-1} \partial_i \partial_j (u_i u_j - \nu \theta \partial_i u_j) =
(-\Delta)^{-1} \partial_i(u_j \partial_j u_i - \nu \partial_j \theta \partial_i u_j) \ .
\label{pressure-1}
\end{equation}
Observe that, as functions of $\theta$, neither the viscosity nor the conductivity
are assumed to be bounded above or below; we only assume that $\theta_0 \geq 1$.
For us, this turns out to be an advantage.
One could imagine that thermal viscosity like in \eqref{thNSE-1} could produce an advanced
regularizing effect for the equation. In regions where the fluid becomes highly turbulent,
the local enstrophy generates more heat via \eqref{Temp-1}. This, in turn, generates extra
viscosity in the trouble-spots. We show that this does indeed happen, and in a way that
can be captured by energy-type estimates weighed against the temperature. In this model,
it is the \emph{possibility} that the viscosity can blow-up with temperature that
prevents turbulent blow-up of the fluid velocity.\\
\\
We remark here that, for most liquids, the viscosity tends to decrease with temperature
\cite{Lautrup}. Although the enhanced Brownian motion effect is still present,
molecules in liquids experience van der Waals forces or hydrogen bonds (the same
forces which cause surface tension), and the added heat disrupts this chemical attraction.
Equations \eqref{thNSE-1} and \eqref{Temp-1} could then be said to model an incompressible,
relatively inert fluid (powder or atmospheric gases). It may be possible to extend
our results to more general profiles $\nu(\theta)$ and $\kappa(\theta)$, so long as
they both go to $\infty$ as $\theta \rightarrow \infty$ (modeling conventional
liquids with an unconventional ``local evaporation'' at high temperatures). This will
be the subject of future work.\\
\\
A nonnegative function $\omega: \TT^3 \rightarrow [0,\infty)$ is said to belong
to $A_2$ if
$$ Q_2(\omega) := \sup_{Q} \left( \frac{1}{|Q|}\int_Q \omega(x) dx \right)
\left(\frac{1}{|Q|} \int_Q \omega(x)^{-1} dx \right) < \infty \ , $$
where the supremum is taken over all cubes $Q \subseteq \TT^3$.
Such functions are called Muckenhoupt weights \cite{HMW}, and play
an important role in weighted Sobolev spaces.\\
\\
\noindent{}With this, we introduce our first main result,
establishing a strong a priori bound on smooth solutions to the above
system, assuming a global bound on $\theta$. We prove this in Section 2.\\

\begin{theorem}[Bound for Thermal Navier-Stokes Eq.]
Let $u$ and $\theta$ be classical solutions to the system \eqref{thNSE-1} and \eqref{Temp-1}
on $\TT^3 \times [0,T]$ with initial data $u_0 \in H^1(\TT^3)$ and $\theta_0 \in L^2(\TT^3)$
with $\theta_0 \geq 1$.
Assume that
\begin{equation}
Q_2 \left( \theta(\cdot, t)^{3/2} \right) < M
\label{muck}
\end{equation}
for some $M>0$ and all $t \in [0,T]$.
Then there exists a constant $L = L(M) > 0$ (see the remark below) such that if
\begin{equation}
K \geq \max(L, 2) \ ,
\label{K-condition}
\end{equation}
then, for another constant $C = C(L) > 0$, we have that
\begin{equation}
\int_{\TT^3} \theta(x,t)^{-1/K} |\nabla u(x,t) |^2 dx \leq
\int_{\TT^3} \theta_0(x)^{-1/K} |\nabla u_0(x)|^2 dx
+ C t \frac{K^8}{\nu^{15}} \left( \frac{1}{2} \| u_0 \|_{L^2}^2 +
\| \theta_0 \|_{L^2}^2 \right)^7
\label{main-1}
\end{equation}
holds for all $t \in [0,T]$. Moreover, the quantities
\begin{equation}
\int_{\TT^3} \theta^{(K-1)/K} |\nabla^2 u|^2 dx \ , \
\int_{\TT^3} \theta^{-(K+1)/K} |\nabla u|^4 dx \ , \
\int_{\TT^3} \theta^{-(K+1)/K} |\nabla u|^2 |\nabla \theta|^2 dx
\label{main-conserved-1}
\end{equation}
are time integrable on $[0,T]$.
\label{thm-main-1}
\end{theorem}
\noindent{}Aside from \eqref{muck}, no additional assumptions of smallness are made on
the solutions $u$ and $\theta$.\\
\\
\textbf{Remark 1:} The constant $L$ emerges largely from bounds on Riesz operators.
We also mention here that $L$ in fact depends on $\inf \theta_0$.
Theorem \ref{thm-main-1} applies equally well to any initial
$u_0$ divergence-free and $\theta_0$ strictly positive, however the dependencies in the
constants become more complicated. Theorem \ref{thm-main-1} could likely also be
extended to include cases where $\theta_0 \geq 0$. Since the temperature essentially
satisfies an advected heat equation with nonnegative forcing, one would expect that
$\theta^2$ becomes strictly positive instantaneously. For this article though, we insist
that $\theta_0 \geq 1$ to simplify the argument.\\
\\
Observe that Theorem \ref{thm-main-1} holds for arbitrarily large conductivity.
This is intuitively
unexpected, since the effect of ``generating extra viscosity in trouble-spots'' is weaker
if the fluid can disperse the ambient temperature more effectively throughout itself;
turbulent areas are cooler when $\kappa$ is larger, and this is seen in the fact that
\eqref{main-1} degenerates as $K \rightarrow \infty$. But the bound still holds, and
the thermal weight even improves with larger $K$ (see Proposition \ref{prop-1}).
This seems to indicate that the build-up of heat in turbulent regions is not indispensable to
the boundedness of the solution. However, taking $\kappa$ to $\infty$ in the thermal
Navier-Stokes system formally leads to the system
\begin{align}
& \partial_t u + u \cdot \nabla u + \nabla p  - \nu \Theta \Delta u = 0
\ , \ \ \ \ \ \ \ \ \ \nabla \cdot u = 0 \nonumber \\
& \Theta(t) = \sqrt{\Theta_0^2 + \frac{1}{2}\| u_0 \|_{L^2}^2
- \frac{1}{2}\| u \|_{L^2}^2} \nonumber
\end{align}
where the effective viscosity is now time-dependent, but constant in space; in particular,
condition \eqref{muck} would be satisfied trivially. Methods of the type
employed in this article seem unlikely to prove similar bounds for the above system.
See, for instance, the comments of \cite{Tao}.\\
\\
\noindent{}To further showcase the effects of thermal viscosity, we also
examine a compressible coupled system of nonlinear fluid-like PDE modeled after the
Burgers equation. The usual multi-dimensional viscous Burgers equation takes the form
$$ \partial_t u + u \cdot \nabla u - \nu \Delta u = 0 \ . $$
Although it lacks compressibility and is vector-valued, the above equation has global
smooth solutions (see for example \cite{Unter}) because it satisfies an $L^\infty$-maximum
principle: taking a dot product of the viscous Burgers equation with $u$ leads to a
drift-diffusion equation for the scalar $|u|^2$, which immediately implies a control
on $\| u \|_{L^\infty}$ for all time. This is sufficient to bootstrap to higher regularity
\cite{Vasseur}.\\
\\
The na\"{i}ve thermal Burgers equation (where we replace $-\nu \Delta u$ by $-\nu
\text{div} (\theta \nabla u)$ and append a transport equation for $\theta^2$ as in
\eqref{Temp-1}) satisfies a similar maximum principle. As such,
we prefer to work on a model for which the base equation is not known to be globally well-posed.
We therefore introduce the ``reduced'' Burgers equation:
$$ \partial_t u +u \cdot \nabla u +\frac{1}{\gamma} (\nabla \cdot u) u -\nu \Delta u = 0 \ , $$
for $2 \leq \gamma < \infty$.
To the best of our knowledge, the above equation has not been shown to be globally
well-posed for smooth initial data on $\TT^3$ (or, for that matter, considered at all).
It no longer satisfies an a priori $L^\infty$ bound since we cannot guarantee any fixed
sign for the divergence at the maximum of $|u|^2$. However, integrating the equation against
$|u|^{\gamma-2} u$ shows that it retains a uniform bound on $\| u \|_{L^\gamma}$ (plus
time-integrability of some other positive quantities coming from the dissipation);
the Burgers equation has been reduced from $L^\infty$ to $L^\gamma$. In the special
case $\gamma = 2$, the a priori estimates are no better than for the Navier-Stokes equations.\\
\\
So, to see the benefit of thermal viscosity in a compressible context (which is closer
to the original Navier-Stokes-Fourier system),
we also examine the initial value problem
\begin{align}
& \partial_t u_j +  u_i \partial_i u_j + \frac{1}{\gamma} u_j \partial_i u_i -
\nu \partial_i (\theta \partial_i u_j) =
0   \label{CTU-F} \\
& \partial_t (\theta^2) + \partial_i ( u_i \theta^2) - \kappa \partial_i (\theta
\partial_i (\theta^2)) = \nu \theta | \nabla u |^2  \label{CTU-T}
\end{align}
on the torus $\TT^3$, with given (sufficiently regular) initial data $u_0$ and $\theta_0$, and
a parameter $\gamma$ as above. We require that $\theta_0(x) \geq 1$ for all $x \in \TT^3$.
Here $\theta^2$ represents the temperature of the system (advected by $u$,
conducted within the medium, and fed by the local dissipation of kinetic energy).\\
\\
Our second (and stronger) main result establishes a similar strong a priori
bound for the exotic thermal system \eqref{CTU-F} and \eqref{CTU-T}, but does not
assume any global bound on $\theta$. We prove this in Section 3.\\

\begin{theorem}[Bound for Thermal Reduced Burgers Eq.]
Assume $\gamma = 2$ or $\gamma \geq 9/4$, and that $\kappa \geq \nu$.
Let $u$ and $\theta$ be classical solutions to the system
\eqref{CTU-F} and \eqref{CTU-T} on $\TT^3 \times [0,T]$ with initial data $u_0 \in H^1(\TT^3)$
and $\theta_0 \in L^2(\TT^3)$ with $\theta_0 \geq 1$.
Then there exist constants $C = C(\nu,\gamma,K,
\| u_0 \|_{L^\gamma}, \| \theta_0 \|_{L^2})>0$ and $M = M(K,\gamma)>0$ such that
\begin{equation}
\int_{\TT^3} \theta(x,t)^{-1/K} | \nabla u (x,t) |^2 dx \leq
\int_{\TT^3} \theta_0(x)^{-1/K} | \nabla u_0(x) |^2 dx + C \left( t^M+1 \right) \ ,
\label{main-general}
\end{equation}
holds for all $t \in [0,T]$. Moreover, the quantities
\begin{equation}
\int_{\TT^3} \theta^{(K-1)/K} |\nabla^2 u|^2 dx \ , \
\int_{\TT^3} \theta^{-(K+1)/K} |\nabla u|^4 dx \ , \
\int_{\TT^3} \theta^{-(K+1)/K} |\nabla u|^2 |\nabla \theta|^2 dx
\label{main-conserved}
\end{equation}
are time integrable on $[0,T]$.
\label{Thm-main}
\end{theorem}
\noindent{}We stress here that the size of the initial data $u_0$ and $\theta_0$ (measured,
for instance, in any Sobolev norm) can be \emph{arbitrarily large} and that the viscosity
$\nu$ can be \emph{arbitrarily small}. The conductivity $\kappa$ can also be arbitrarily small,
so long as it remains bigger than $\nu$.\\
\\
\noindent{}The constants above all depend on $\inf \theta_0$ in nontrivial ways.
The polynomial growth in \eqref{main-general} comes from a weakened thermodynamic
lower bound on the temperature (see Proposition \ref{Prop-initial}).
The condition on the relative sizes of $\kappa$ and $\nu$ is slightly weaker;
we only need $K \geq 2$. This condition emerges in an
interesting way (that is also present in the proof of Theorem \ref{thm-main-1}). When we try
to do the
weighted enstrophy estimate, the diffusion itself produces terms which are roughly like
$$ \int \theta^{-\frac{K+1}{K}} |\nabla u|^2 |\nabla \theta|^2 dx \ . $$
In fact, one such term (named $J_1$ in both proofs below) has a definite (and detrimental)
sign in the resulting inequality. So special cancellations should not be expected.\\
\\
We also remark
that the actual lower bound for $K$ implied by our methods is a bit lower than $2$, but then
the dependencies in the constants become more complicated. Theorem \ref{Thm-main} would
remain true with
$$ K > \frac{3+\sqrt{17}}{4} \approx 1.78 \ , $$
and condition \eqref{K-condition} could be weakened to
$$ K > \max(L, (3+\sqrt{17})/4) \ . $$
Essentially, we only need the term labeled $U_2$ in inequalities \eqref{ineq-1} and
\eqref{final-ineq} below to have a positive coefficient. A slightly different approach
(see Remark 2) improves this lower bound further, but it seems that the large Prandtl number
regime ($K \approx 1$) will require some fundamentally different approaches.\\
\\
Although the proof of Theorem \ref{Thm-main} is similar to the proof for
Theorem \ref{thm-main-1}, there are two aspects
that we believe make it worth presenting. Firstly,
we were rather surprised that compressibility did not hinder the main argument,
since it is fundamentally an energy-type estimate. Indeed, the pressure term in
the incompressible model is what forces the extra assumption \eqref{muck}
on the temperature. Yet, for this model, the result holds with no additional assumptions.
Secondly, the initial thermodynamic a priori estimates
are harder to prove and less clear-cut than for the thermal Navier-Stokes system; compare
Proposition \ref{prop-1} with Proposition \ref{Prop-initial}. Compressibility allows
for a ``refrigeration'' effect,
where an expanding gas becomes colder. This is seen explicitly in \eqref{thet*t}.
Nevertheless, we are still able to prove a (time-decaying) lower bound on the temperature
and an upper bound on the total heat $\| \theta \|_{L^2}^2$ for the thermal reduced
Burgers equation (both of which are crucial). Beyond Proposition \ref{Prop-initial},
the rest of the proof of Theorem \ref{Thm-main} is no harder than for the
incompressible case.\\
\\
Bootstrapping the above a priori bounds to higher regularity can be done through classical
techniques \cite{Vasseur}. This is expected, since the new conserved quantities easily
satisfy the Ladyzhenskaya-Prodi-Serrin condition (\cite{Ser}, \cite{Lady}), ignoring for the
moment that the dissipation has potentially rough coefficients in the form of $\theta$.
A more modern (or direct) approach would be along the lines of \cite{S-V}, where
H\"{o}lder continuity is propagated through the use of Campanato's characterization of
H\"{o}lder spaces \cite{Camp}: for $\phi \in C_c^\infty(\RR^n)$ nonnegative, radially
symmetric, supported in $B_1(0)$, with $\int \phi(y) dy = 1$ and $f \in L^2(\RR^n)$, \emph{if}
there are $A>0$ and $\alpha \in (0,1)$ such that for all $r>0$ and $x \in \RR^n$,
$$ \int_{B_1(0)} \left| f(x+ry) - \int_{B_1(0)} f(x+rz) \phi(z) dz \right|^2
\phi(y) dy \leq A^2 r^{2 \alpha} \ , $$
\emph{then} $f$ has a H\"{o}lder continuous representative in $C^\alpha(\RR^n)$ (and
$[f]_{C^\alpha} \lesssim A$). The conserved
quantities of \eqref{main-conserved-1} and \eqref{main-conserved} should suffice to repeat the
proof of \cite{S-V} for $u$. Note that, since
$$ \int \theta^{(K-1)/K} |\nabla^2 u|^2 dx $$
is time-integrable, then it is in particular finite for a sequence of times that converge to
$0$. By Sobolev embedding ($H^2(\TT^3) \subset C^{1/2}(\TT^3)$), the propagation argument of
\cite{S-V} would imply our solutions $u$ become H\"{o}lder continuous instantly, even with
$u_0 \in H^1$.\\
\\
Although we only present a strong a priori bound for classical solutions, the global
well-posedness of the Cauchy problems \eqref{thNSE-1}-\eqref{Temp-1} and
\eqref{CTU-F}-\eqref{CTU-T} in Sobolev space should follow from established techniques.
The existence of weak solutions could be proven in a similar manner as for the
Navier-Stokes-Fourier system (see \cite{D-F} and \cite{F-N}). From there, one could prove
a local well-posedness result (as in \cite{M-N} for the Navier-Stokes-Fourier system) to
obtain classical solutions on a
short time interval that depends sensitively on the initial data. Theorems \ref{thm-main-1}
and \ref{Thm-main} would then extend this solution to a global one. Alternatively,
one could mollify the weak solutions and use the results of this paper to show a convergence
of mollified solutions to a unique strong solution.
The complete proof of global well-posedness is an interesting topic and will be the subject
of a forthcoming paper.\\

\section{Thermal NSE}
\noindent{}We first prove some essential thermodynamic properties of the system \eqref{thNSE-1}
and \eqref{Temp-1}.
\begin{proposition}[Initial bounds on energy and minimal temperature]
For $u$ and $\theta$ smooth solutions to \eqref{thNSE-1} and \eqref{Temp-1} on the time
interval $[0,T]$, with initial data $u_0$ and $\theta_0 \geq 1$ respectively, we have
$$ \inf_{\TT^3} \theta( \cdot ,t) \geq 1 \ , $$
for all $t \in [0,T]$ and
$$ \sup_{t \in [0,T]} \left( \frac{1}{2} \| u(\cdot, t) \|_{L^2(\TT^3)}^2 +
\| \theta(\cdot, t) \|_{L^2(\TT^3)}^2 \right) \leq
\frac{1}{2} \| u_0 \|_{L^2(\TT^3)}^2 + \| \theta_0 \|_{L^2(\TT^3)}^2 \ . $$
\label{prop-1}
\end{proposition}
\noindent{\sc Proof.} The minimum principle for $\theta$ is intuitive, but we will demonstrate
it rigorously. Take $\hat{T} > 0$ and $\epsilon > 0$ to be determined later. Define
$v(x,t) = \theta(x,t) + \epsilon t$. Then, by \eqref{Thet-1},
\begin{equation}
\partial_t v + u \cdot \nabla v - \kappa \theta \Delta v
 - 2 \kappa |\nabla v|^2 = \frac{\nu}{2} |\nabla u|^2 + \epsilon \ .
 \label{thet+eps}
 \end{equation}
Since $v$ is a smooth function and $\TT^3 \times [0,\hat{T}]$ is compact, it achieves a global
minimum; say this happens at $(\hat{x},\hat{t})$. If $\hat{t} > 0$, then we know that
$$ \partial_t v(\hat{x},\hat{t}) \leq 0 \ \ , \ \ \Delta v(\hat{x},\hat{t}) \geq 0
\ \ \text{, and} \ \ \nabla v(\hat{x},\hat{t}) = 0 \ . $$
However, so long as $\theta(x,t) \geq 0$, \eqref{thet+eps} shows that
$$ \partial_t v(\hat{x},\hat{t}) \geq \epsilon \ . $$
Let $\tilde{T}$ be the maximum time for which $\theta \geq 0$ on
$\TT^3 \times [0, \tilde{T}]$.
This is positive by continuity and the assumption that $\theta_0(x) \geq 1$.
Note that this is also a time interval on which \eqref{Thet-1} is valid.
Choose $\hat{T} = \min (\tilde{T}, T)$. Then the minimum of $v$ on
$\TT^3 \times [0, \hat{T}]$ must happen at time $\hat{t}=0$. But, since $\epsilon$ was
arbitrary, this means
$$ \inf_{\TT^3 \times [0,\hat{T}]} \theta(x,t) = \inf_{\TT^3} \theta_0(x) \geq 1 \ . $$
In particular, this shows that $\tilde{T} > \hat{T}$ (again by continuity) and that $\hat{T}=T$,
proving the first part of Proposition \ref{prop-1} (and that \eqref{Thet-1} is an effective
substitute for \eqref{Temp-1}).\\
\\
Integrating \eqref{thNSE-1} against $u_j$ and integrating \eqref{Temp-1} by itself yields
$$ \frac{\partial_t}{2} \| u \|_{L^2}^2 + \nu \int \theta |\nabla u|^2 dx = 0 \ \ \ \ \ 
\text{ and } \ \ \ \ \ \partial_t \| \theta \|_{L^2}^2 = \nu \int \theta |\nabla u|^2 dx \ . $$
Adding the two shows that
$$ \partial_t \left( \frac{1}{2} \| u \|_{L^2}^2 + \| \theta \|_{L^2}^2 \right) = 0 \ , $$
which completes the proof of the Proposition.
\indent \hfill $\bullet$\\
\\
It is also true that we get an a priori control on the ``thermal enstrophy''
$$ \nu \int_0^T \int \theta |\nabla u|^2 dx dt \leq \frac{1}{2} \| u_0 \|_{L^2}^2 \ , $$
as with the regular Navier-Stokes equations. Interestingly, the proof of the main theorem
will not require such a bound, though it establishes much stronger conserved quantities.\\
\\
\noindent{}{\sc Proof of Theorem 1.1.}
Let $f: \RR_+ \rightarrow \RR_+$ (to be specified later) be our thermal weight
function. We proceed by integrating \eqref{thNSE-1} against
$-\text{div}(f(\theta) \nabla u)$ to get:
$$  \underbrace{\int f \frac{\partial_t}{2} |\nabla u|^2 dx}_{I_0}
\underbrace{- \int u_i \partial_i u_j \partial_k (f \partial_k u_j) dx}_{I_A}
\underbrace{- \int \partial_j p \partial_k (f \partial_k u_j) dx}_{I_P}
\underbrace{+ \nu \int \partial_i (\theta \partial_i u_j)
\partial_k (f \partial_k u_j) dx}_{I_D} = 0  $$
\noindent{}The advection term is straightforward:
\begin{equation}
I_A = \int \partial_k u_i \partial_i u_j f \partial_k u_j dx +
\int u_i \partial_i \partial_k u_j f \partial_k u_j dx \geq
-\int f |\nabla u|^3 dx - \frac{1}{2} \int f' |\nabla u|^2 u_i \partial_i \theta dx
\label{IA-1}
\end{equation}
Recalling \eqref{pressure-1}, we integrate by parts twice in $I_P$ to get the formula
\begin{equation}
I_P = -\int R_l R_k [u_m \partial_m u_l - \nu \partial_m \theta \partial_l u_m]
f' \partial_j \theta \partial_k u_j dx \ ,   \label{IP-1}
\end{equation}
where $R_j$ is the Riesz operator $\partial_j (-\Delta)^{-1/2}$ on the torus, with
Fourier symbol $i \frac{k_j}{|k|}$.\\
\\
Lastly, $I_D$ will provide three positive terms (and one cancellation) for our estimates.
\begin{align}
I_D &= \nu \int \partial_k(\theta \partial_i u_j) \partial_i (f \partial_k u_j) dx =
\nu \int (\partial_k \theta \partial_i u_j +
\theta \partial_k \partial_i u_j)(f' \partial_i \theta \partial_k u_j +
f \partial_i \partial_k u_j) dx  \nonumber \\
&= \nu \int \left( f' |\nabla \theta \cdot \nabla u|^2 +
f \partial_k \theta \partial_i u_j \partial_i \partial_k u_j +
\theta f' \partial_i \theta \partial_k u_j \partial_k \partial_i u_j +
\theta f |\nabla^2 u|^2 \right) dx \nonumber \\
&= \underbrace{\nu \int f' |\nabla \theta \cdot \nabla u|^2  dx}_{J_1} +
\underbrace{\nu \int (f + \theta f') \partial_k \theta \partial_i u_j
\partial_i \partial_k u_j dx}_{J_2} + \underbrace{\nu \int \theta f |\nabla^2 u|^2 dx}_{J_3}
\nonumber
\end{align}
The last equality was obtained by re-indexing. Observe that $J_3$ is a positive term.
Continuing, we integrate $J_2$ by parts to move the $k$-derivative away from
$\partial_i u_j$ and onto $\partial_k \theta$:
\begin{align}
J_2 &= -\frac{\nu}{2} \int (f+\theta f') \partial_k \partial_k \theta |\nabla u|^2 dx -
\frac{\nu}{2} \int \partial_k (f+\theta f') \partial_k \theta |\nabla u |^2 dx  \nonumber \\
&= \frac{\nu}{2 \kappa} \int \frac{f+\theta f'}{\theta} (-\kappa \theta \Delta \theta)
|\nabla u|^2 dx -
\underbrace{\frac{\nu}{2} \int (2f' + \theta f'') |\nabla \theta|^2 |\nabla u|^2 dx}_{K_0}
\nonumber \\
&= \frac{1}{2} \int \left( \frac{\nu}{\kappa}\frac{f+\theta f'}{\theta} \right)
\left(\underbrace{\frac{\nu}{2}|\nabla u|^2}_{K_1}
\underbrace{+ 2 \kappa | \nabla \theta|^2}_{K_2} \underbrace{- u \cdot \nabla \theta}_{K_3}
\underbrace{- \partial_t \theta}_{K_4} \right) |\nabla u|^2 dx - K_0
\nonumber
\end{align}
The last equality comes from direct substitution with \eqref{Thet-1}. We would like to
merge the integral involving the time derivative of $\theta$ (that is, the one coming
from the $K_4$ term) with $I_0$ to pull the time derivative entirely outside the integral
(to recover the weighted enstrophy norm of \eqref{main-1}).
In order to do this, we must clearly have
\begin{equation}
f'(\theta) = -\frac{\nu}{\kappa} \frac{f(\theta) + \theta f'(\theta)}{\theta}
\ \ \ \ \text{ or } \ \ \ \ f(\theta) = K \theta^{-\frac{1}{K}}
\ \ \ \text{ and } \ \ \ f'(\theta) = -\theta^{-\frac{K+1}{K}} \ . \nonumber
\end{equation}
\textbf{Remark 2:} The more obvious approach for computing energy estimates with weights
would be to start with $\int f(\theta) |\nabla u|^2 dx$, take a time derivative,
then use equations \eqref{thNSE-1} and \eqref{Thet-1} to obtain the inequality. This approach
does not, at first, impose any condition on the weight $f$. The method above is basically the
same, but we demanded that the term arising from $-\kappa \theta \Delta \theta$ in
\eqref{Thet-1} exactly cancels with a term (containing $\Delta \theta$) that arises from
integrating the dissipation of \eqref{thNSE-1} by parts; this fixes $f$ as above.
We mention that, using the conventional approach outlined in this remark, we could not
produce a better final inequality or even a different class of admissible weights.\\
\\
Having established $f(\theta)$ explicitly, we now make some observations. As mentioned before,
$$ I_0 + \frac{1}{2} \int \left( \frac{\nu}{\kappa} \frac{f+\theta f'}{\theta} \right) K_4
|\nabla u|^2 dx = \frac{\partial_t}{2} \int K \theta^{-\frac{1}{K}} | \nabla u|^2 dx. $$
Moreover, the second term in \eqref{IA-1} from $I_A$ exactly cancels with the integral
arising from $K_3$.\\
\\
\noindent{}Also, $K_0$ and the integral arising from $K_2$ combine into
\begin{equation}
\nu \left( K - \frac{1}{2} - \frac{1}{2K} \right) \int \theta^{-\frac{K+1}{K}}
|\nabla \theta|^2 |\nabla u|^2 dx \label{KK-1}
\end{equation}
which is also positive.\\
\\
\noindent{}$J_1$ is evidently negative. Observe, however, that it is controlled by
\eqref{KK-1} (since $\sum_{i,j} a_i b_i a_j b_j \leq \sum_{i,j} a_i^2 b_j^2$), provided
$K - 1/(2K) > 3/2$; which certainly holds if $K \geq 2$. This is the first time we see the
methodology's dependence on $K$.\\
\\
\noindent{}The crucial benefit from this substitution is seen in the integral term arising
from $K_1$. This gives us a gradient of $u$ to the \emph{fourth} power with a positive sign
(weighed by a negative power of $\theta$ that improves as $K$ increases). It is here that
we see the heat production of \eqref{Thet-1} generate an improvement in the
dissipation of \eqref{thNSE-1}. Combined with
\eqref{KK-1}, these two positive terms allow us to close the estimate in a novel way.\\
\\
\textbf{Remark 3:} We mention here that the above method requires that the viscosity
and conductivity both grow with temperature. If, in place of \eqref{thNSE-1}
and \eqref{Temp-1}, we consider the system
\begin{align}
& \partial_t u + u \cdot \nabla u + \nabla p - \nu \text{div}(w(\theta) \nabla u) = 0
\nonumber \\
& \partial_t (\theta^2) + u \cdot \nabla (\theta^2) -
\kappa \text{div}(w(\theta) \nabla (\theta^2)) =
\nu w(\theta) | \nabla u |^2  \nonumber
\end{align}
for $w$ a smooth positive weight, then we could proceed in a similar fashion up until
the analysis of the term $J_2$. In order for the cancellation to occur properly, we would
need
$$ \frac{\nu}{\kappa} \frac{f w' + f' w}{w} = -f' \ , $$
so that (up to a multiplicative constant)
$$ f(\theta) = K w(\theta)^{-1/K} \ . $$
If $f'$ were not strictly negative, we would no longer have a gain from the $K_1$ and
$K_2$ terms; indeed, those terms would have a bad sign if $f$ were increasing in
$\theta$. Since
$$ f'(\theta) = - w(\theta)^{-(K+1)/K} w'(\theta) \ , $$
we see that the above improvement can only happen if $w' > 0$.\\
\\
\noindent{}Collecting what terms remain, we obtain the inequality:
\begin{align}
\frac{\partial_t}{2} \int K \theta^{-\frac{1}{K}} & |\nabla u|^2 dx +
\underbrace{\nu K \int \theta^{\frac{K-1}{K}} |\nabla^2 u|^2 dx}_{U_1} +
\underbrace{\nu \frac{2K^2-3K-1}{2K} \int \theta^{-\frac{K+1}{K}} |\nabla \theta|^2
|\nabla u|^2 dx}_{U_2} \nonumber \\
&+ \underbrace{\frac{\nu}{4} \int \theta^{-\frac{K+1}{K}} |\nabla u|^4 dx}_{U_3}
\leq K \int \theta^{-\frac{1}{K}} |\nabla u|^3 dx - I_P
\label{ineq-1}
\end{align}
We use H\"{o}lder's inequality and Sobolev embedding in three dimensions to control
the cubic term on the right hand side.
\begin{align}
K \int \theta^{-\frac{1}{K}} |\nabla u|^3 dx &\leq
K \left(\int |\nabla u|^6 dx \right)^{\frac{1}{9}} \left( \int |\nabla u|^4
\theta^{-\frac{K+1}{K}} dx \right)^{\frac{7}{12}} \left( \int \theta^2 dx
\right)^{\frac{7K-5}{24K}} | \TT^3 |^{\frac{7K-3}{24K}} \nonumber \\
& \leq C K \left( \int | \nabla^2 u|^2 dx \right)^{\frac{1}{3}}
\left( \int |\nabla u|^4 \theta^{-\frac{K+1}{K}} dx \right)^{\frac{7}{12}}
\left( \int \theta^2 dx
\right)^{\frac{7K-5}{24K}} \nonumber
\end{align}
where $C$ incorporates the Sobolev constant and the size of the torus; taking $C$ large
enough also makes it independent of $K$.
Using the minimum principle for $\theta$ and Young's inequality with exponents $3$,
$\frac{12}{7}$, and $12$, we see that
\begin{equation}
K \int \theta^{-\frac{1}{K}} |\nabla u|^3 dx \leq
\frac{1}{3} U_1 + \frac{1}{3} U_3
+ C \frac{K^8}{\nu^{11}} \| \theta \|_{L^2}^{\frac{7K-5}{K}}
\label{cubic-term-1}
\end{equation}
for a different constant $C$ independent of $K$ and $\nu$. Also recall that
$$ \| \theta \|_{L^2}^2 \leq \| \theta_0 \|_{L^2}^2 + \frac{1}{2} \| u_0 \|_{L^2}^2 \ . $$
Thus, the cubic term on the right hand side of \eqref{ineq-1} is indeed dominated
by the positive terms and the conserved energy.\\
\\
\noindent{}To estimate the term arising from the pressure, we quote a result concerning
optimal bounds for the operator norm of Riesz transforms
on weighted Lebesgue spaces. For a weight $\omega: \TT^3 \rightarrow [0,\infty)$, the
weighted $L^2$-norm of a function $\phi$ is given by
$$ \| \phi \|_{L^2(\omega)} = \left( \int_{\TT^3} \phi(x)^2 \omega(x) dx \right)^{1/2} \ .$$
Theorem 2.1 of \cite{Pet} states that there exists a constant $c$ so that for all weights
$\omega \in A_2$, the Riesz transforms $R_k: L^2(\omega) \rightarrow L^2(\omega)$ have
operator norm $\| R_k \| \leq c Q_2(\omega)$.\\
\\
By assumption \eqref{muck}, $Q_2(\theta(\cdot, t)^{3/2}) < M$ for all $t$. Since $K \geq 2$,
we know that $(K+1)/K \leq 3/2$. Furthermore, $A_2$ is closed under convex interpolation:
if $\omega$, $\sigma \in A_2$, then for all $\alpha \in (0,1)$, we have that
$\omega^\alpha \sigma^{1-\alpha} \in A_2$
and
$$ Q_2(\omega^\alpha \sigma^{1-\alpha}) \leq Q_2(\omega)^\alpha Q_2(\sigma)^{1-\alpha} \ , $$
by H\"{o}lder's inequality. Therefore, $Q_2(\theta(\cdot, t)^{(K+1)/K}) < M$
for all $t \in [0,T]$ (interpolating between $\theta^{3/2}$ and $1$), and we conclude that
\begin{equation}
\| R_k R_l g \|_{L^2(\theta^{(K+1)/K})} \leq c M^2 \| g \|_{L^2(\theta^{(K+1)/K})} \ ,
\label{weig-riesz}
\end{equation}
for all smooth functions $g$ and a constant $c$ independent of all parameters.\\
\\
Now we estimate $I_P$. From \eqref{IP-1} and the explicit form of $f$, we
have that
\begin{align}
|I_p| &= \left| \int R_k R_l (u_m \partial_m u_l - \nu \partial_m \theta \partial_l u_m)
\theta^{-\frac{K+1}{K}} \partial_j \theta \partial_k u_j dx \right| \nonumber \\
&\leq \frac{C}{M\nu}\int\left| R_k R_l(u_m\partial_m u_l)\right|^2 \theta^{-\frac{K+1}{K}}dx
+ M\nu \int |\nabla \theta|^2 |\nabla u|^2 \theta^{-\frac{K+1}{K}}dx \nonumber \\
& \ \ \ + \frac{\nu}{M} \int\left| R_k R_l(\partial_m \theta \partial_l u_m) \right|^2
\theta^{-\frac{K+1}{K}}dx + \bar{C}M\nu \int |\nabla \theta|^2 |\nabla u|^2
\theta^{-\frac{K+1}{K}}dx \nonumber \\
&\leq \frac{CM}{\nu} \int |u|^2 |\nabla u|^2 \theta^{-\frac{K+1}{K}}dx+
\bar{C} M \nu \int |\nabla \theta|^2 |\nabla u|^2 \theta^{-\frac{K+1}{K}}dx \ , \nonumber
\end{align}
where $\bar{C}$ depends on $c$ from \eqref{weig-riesz} and the dimension. Therefore,
\begin{equation}
|I_P| \leq \frac{C}{\nu} \int \theta^{-\frac{K+1}{K}} |u|^2 |\nabla u|^2 dx
+ \frac{2}{3} U_2 \ ,
\label{pressure-semi-1}
\end{equation}
assuming \eqref{K-condition}. It is here that we determine the constant in
\eqref{K-condition} (and see that it is essentially proportional to $M$). The final constant
$C$ in the inequality above incorporates the factor of $M$ from the previous estimate; the
remaining bounds will not impose any further restrictions on $K$.\\
\\
\textbf{Remark 3:} There is a somewhat different approach to bounding $I_P$, involving
commutator estimates. We quote a second (classical) result by
\cite{CRG}: for $R_k$ a Riesz operator and $\phi , g \in C^\infty_0(\RR^n)$, if
$$ [R_k, \phi] g := R_k (\phi g) - \phi R_k g \ , $$
then
$$ \| [R_k, \phi] g \|_{L^p} \lesssim [\phi]_{BMO} \| g \|_{L^p} \ , $$
where the constant implied by the inequality only depends on $p$ and the
dimension $n$.
Theorem I of \cite{CRG} covers more general Calder\'{o}n-Zygmund singular
integral operators (such as the composite $R_k R_l$), but for our purposes we quote the
above bound and observe that
\begin{equation}
[R_k R_l , \phi](g) = R_k \left( [R_l,\phi](g) \right) +
[R_k,\phi] \left( R_l g \right) \ . \nonumber
\end{equation}
By boundedness of Riesz transforms on $L^p$ for $1<p<\infty$ (see for instance \cite{Stein}),
we see that
\begin{equation}
\| [ R_k R_l , \phi](g) \|_{L^2(\TT^3)} \leq c [\phi]_{BMO} \| g \|_{L^2(\TT^3)} \ ,
\label{comm-estimate}
\end{equation}
with $c$ independent of $\phi$ or $g$.\\
\\
We could then write
\begin{align}
I_P &= \int R_k R_l \left( \theta^{-\frac{K+1}{2K}} (u_m \partial_m u_l - \nu \partial_m \theta
\partial_l u_m) \right) \theta^{-\frac{K+1}{2K}} \partial_j \theta \partial_k u_j dx \nonumber \\
& \ \ \ +\int [R_k R_l, \theta^{\frac{K+1}{2K}}] \left( \theta^{-\frac{K+1}{2K}}
(u_m \partial_m u_l - \nu \partial_m \theta \partial_l u_m) \right) \theta^{-\frac{K+1}{K}}
\partial_j \theta \partial_k u_j dx \ , \nonumber
\end{align}
and conclude an identical bound to \eqref{pressure-semi-1}, if we assumed that
$$ [\theta(\cdot,t)^{(K+1)/(2K)}]_{BMO} < M $$
for all $t \in [0,T]$. We do not mention this possibility in Theorem \ref{thm-main-1}, in part
because it is inelegant ($K$ must be larger than $L(M) \approx \sup [\theta^{(K+1)/(2K)}]_{BMO}$), but
largely because it would be a stronger type of assumption than \eqref{muck}. A result of
\cite{JN} states that any positive BMO function whose reciprocal is also BMO actually belongs
to $A_p$ for all $p \in (1, \infty]$; recall that $\theta^{-(K+1)/(2K)} \in L^\infty \subset
BMO$.\\
\\
To complete the bound, we estimate
$$ \int \theta^{-\frac{K+1}{K}} |u|^2 |\nabla u|^2 dx \leq
\left( \int |u|^4 dx \right)^{1/2} \left( \int \theta^{-\frac{K+1}{K}}
|\nabla u|^4 dx \right)^{1/2} \| \theta^{-\frac{K+1}{2K}} \|_{L^\infty} \ . $$
Agmon's inequality (Lemma 13.2 of \cite{Agmon}) on a compact three-dimensional domain says that
$$ \| u \|_{L^\infty} \leq C \| u \|_{L^2}^{1/4} \| u \|_{H^2}^{3/4} \ . $$
Thus, by H\"{o}lder's inequality and monotonicity of kinetic energy,
$$ \| u \|_{L^4}^4 \leq \| u \|_{L^2}^2 \| u \|_{L^\infty}^2 \leq C \| u_0 \|_{L^2}^{5/2}
\left( \int \theta^{\frac{K-1}{K}} |\nabla^2 u|^2 dx + \| u_0 \|_{L^2}^2 \right)^{3/4} \ , $$
for yet another constant $C$ independent of all parameters. Therefore
\begin{align}
\frac{C}{\nu} \int \theta^{-\frac{K+1}{K}} |u|^2 |\nabla u|^2 dx & \leq
\frac{\nu}{12} \int \theta^{-\frac{K+1}{K}} |\nabla u|^4 dx +
\frac{C}{\nu^3} \| u_0 \|^{5/2}_{L^2} \left( \int \theta^{\frac{K-1}{K}} |\nabla^2 u|^2 dx
+ \| u_0 \|^2_{L^2} \right)^{3/4} \nonumber \\
& \leq \frac{1}{3} U_3 + \frac{1}{3} U_1 + K \nu \| u_0 \|_{L^2}^2
+ \frac{C}{K^3 \nu^{15}} \| u_0 \|_{L^2}^{10} \ ,
\label{pressure-done-1}
\end{align}
where $C$ has become larger but is still independent of $K$, $\nu$, or the initial data.\\
\\
Combining all three estimates (\eqref{cubic-term-1}, \eqref{pressure-semi-1}, and
\eqref{pressure-done-1}), inequality \eqref{ineq-1} becomes
$$ \frac{\partial_t}{2} \int K \theta^{-\frac{1}{K}} |\nabla u|^2 dx + \frac{1}{3}
(U_1 + U_2 + U_3) \leq \nu K \| u_0 \|_{L^2}^2 +
C \frac{K^8}{\nu^{11}} (\| \theta_0 \|_{L^2}^2 +
\frac{1}{2}\| u_0 \|_{L^2}^2)^{\frac{7K-5}{2K}} +
\frac{C}{K^3 \nu^{15}} \| u_0 \|^{10}_{L^2} \ . $$
Integrating this inequality yields \eqref{main-1} and \eqref{main-conserved-1}.\\
\indent \hfill $\bullet$
\\
\section{Thermal Reduced Burgers}
\noindent{}We begin by proving some thermodynamic estimates that are similar to
Proposition \ref{prop-1}.
The lack of incompressibility, however, will make the a priori bounds more complicated and
will also reduce the strength of the minimum principle for the temperature.\\
\begin{proposition}[Initial bounds on velocity and temperature]
Let $\gamma \geq 2$. For $u$ and $\theta$ smooth solutions to \eqref{CTU-F} and
\eqref{CTU-T} on the time interval $[0,T]$, with initial data $u_0$ and $\theta_0 \geq 1$
respectively, we have
\begin{equation}
\theta(x,t) \geq \frac{1}{3t/(8\nu) + 1} \ ,
\label{CTU-lower}
\end{equation}
for all $(x,t) \in \TT^3 \times [0,T]$ and
\begin{align}
\sup_{t \in [0,T]} & \frac{1}{\gamma} \| u( \cdot, t) \|_{L^\gamma(\TT^3)}^\gamma +
\nu \int_0^T \int_{\TT^3} \theta(x,t) |\nabla u (x,t)|^2 |u(x,t)|^{\gamma-2} dx dt \nonumber \\
&+ \frac{4 \nu (\gamma-2)}{\gamma^2} \int_0^T \int_{\TT^3} \theta(x,t)
| \nabla |u(x,t)|^{\gamma / 2}|^2 dx dt
\leq \frac{1}{\gamma} \| u_0 \|_{L^\gamma(\TT^3)}^\gamma \ .
\label{CTU-gamma}
\end{align}
Moreover, if $\gamma \geq 9/4$ or $\gamma = 2$, we also have
\begin{equation}
\sup_{t \in [0,T]} \left( \frac{1}{2} \| u (\cdot, t) \|_{L^2}^2 + \| \theta(\cdot, t)
\|_{L^2}^2 \right) \leq \frac{1}{2} \| u_0 \|_{L^2}^2 + \| \theta_0 \|_{L^2}^2 +
C(T^2+1) \left( \| u_0 \|_{L^\gamma}^{3\gamma-3} + 1 \right) \ .
\label{CTU-upper}
\end{equation}
\label{Prop-initial}
\end{proposition}
\noindent{\sc Proof.} We begin with \eqref{CTU-lower}, providing
a ``minimum local strength'' for the thermal viscosity. As with the incompressible case,
it will be convenient from now on to work with $\theta$ instead of the actual temperature;
\eqref{CTU-T} then reduces to
\begin{equation}
\partial_t \theta + \frac{1}{2} \text{div}(u) \theta + u_i \partial_i \theta -
\kappa \theta \Delta \theta = 2 \kappa | \nabla \theta |^2 +
\frac{\nu}{2} | \nabla u |^2 \ , \label{CTU-th}
\end{equation}
on any time interval where $\theta \geq 0$.
Let the maximum such time interval be $[0,\tilde{T}]$.
As in the proof of Proposition \ref{prop-1}, $\tilde{T} > 0$ by continuity
and the assumption that $\theta_0 \geq 1$. Define $\hat{T} = \min (T, \tilde{T})$ and define
$$ v(x,t) = \left( \frac{3 t}{8\nu} + 1 \right) \theta(x,t) \ . $$
Then, on $\TT^3 \times [0,\hat{T}]$, $v$ satisfies the equation
\begin{equation}
\partial_t v + u \cdot \nabla v - \kappa \theta \Delta v - \frac{2 \kappa}{3t/(8\nu)+1}
|\nabla v|^2 = \frac{1}{2} \left( \frac{3t}{8\nu}+1 \right) \left( \nu |\nabla u|^2 -
\text{div}(u) \theta \right) + \frac{3}{8\nu} \theta \ .
\label{thet*t}
\end{equation}
Since $v$ is smooth and $\TT^3 \times [0,\hat{T}]$ is compact, $v$ achieves its global
minimum; say at $(\hat{x},\hat{t})$. If $\hat{t} > 0$, then
$$ \partial_t v(\hat{x},\hat{t}) \leq 0 \ \ , \ \ \Delta v(\hat{x},\hat{t}) \geq 0 \ \
\text{, and} \ \ \nabla v (\hat{x},\hat{t}) = 0 \ . $$
By Young's inequality,
$$ \left| \text{div}(u) \theta \right| \leq \frac{\nu}{3} (\text{div}(u))^2 +
\frac{3}{4 \nu} \theta^2 \leq \nu |\nabla u|^2 + \frac{3}{4 \nu} \theta^2 \ . $$
We then have from \eqref{thet*t} that
$$ \partial_t v(\hat{x},\hat{t}) \geq \left( \frac{3t}{8\nu} + 1 \right)
\left( -\frac{3 \theta^2}{8 \nu} \right) + \frac{3 \theta}{8\nu} =
\frac{v(\hat{x},\hat{t})-v(\hat{x},\hat{t})^2}{t + 8\nu / 3} \ . $$
Assume that $v(\hat{x},\hat{t}) < 1$. Then
$\partial_t v(\hat{x},\hat{t})$ is strictly positive, and this can only happen if
$\hat{t}=0$. But then
$$ \inf_{\TT^3 \times [0,\hat{T}]} v(x,t) = \inf_{\TT^3} v(x, 0) =
\inf_{\TT^3} \theta_0(x) \geq 1 \ , $$
which is a contradiction. Therefore $v(x,t) \geq 1$ on $\TT^3 \times [0,\hat{T}]$.
But this shows that $\tilde{T} > \hat{T}$, so that $\hat{T} = T$, and
that \eqref{CTU-lower} holds.\\
\\
To prove \eqref{CTU-gamma}, we multiply \eqref{CTU-F} by $u_j |u|^{\gamma - 2}$, sum in $j$,
and integrate in space. Integration by parts shows that the nonlinearities (in $u$) cancel
each other, leaving
$$ \frac{\partial_t}{\gamma} \| u \|_{L^\gamma}^\gamma +
\nu \int \theta | \nabla u |^2 |u|^{\gamma - 2} dx
+ \frac{4 \nu (\gamma - 2)}{\gamma^2} \int \theta | \nabla |u|^{\gamma / 2} |^2 dx = 0 \ . $$
Recall that $\theta$ stays positive. Integrating this in time from $0$ to $\bar{T}$, then
taking a supremum in $\bar{T}$ over $[0,T]$ yields \eqref{CTU-gamma}.\\
\\
Lastly, our main estimates will require some manner of upper bound on $\theta$. It will
suffice to have a bound on the ``total heat'' in the system, given by $\| \theta \|_{L^2}^2$.
To do this, we multiply \eqref{CTU-F} by $u_j$, sum in $j$, and integrate in space
(as if $\gamma$ were $2$). This gives
$$ \frac{\partial_t}{2} \| u \|_{L^2}^2 + \nu \int \theta | \nabla u |^2 dx =
\left( \frac{1}{2} - \frac{1}{\gamma} \right) \int |u|^2 \text{div}(u) dx \ . $$
If $\gamma = 2$, we integrate \eqref{CTU-T} in space and add it to the above to get
true conservation of kinetic and thermal energy. Then \eqref{CTU-upper} takes on the simpler form
$$ \frac{1}{2} \| u (\cdot, t) \|_{L^2}^2 + \| \theta (\cdot, t) \|_{L^2}^2 =
\frac{1}{2} \| u_0 \|_{L^2}^2 + \| \theta_0 \|_{L^2}^2 \ . $$
If $\gamma \neq 2$, we integrate \eqref{CTU-T} in space and add half of it to the above,
yielding
$$\frac{\partial_t}{2} \left( \| u \|_{L^2}^2 + \| \theta \|_{L^2}^2 \right)
+ \frac{\nu}{2} \int \theta |\nabla u|^2 dx =
\frac{\gamma-2}{2 \gamma} \int |u|^2 \text{div}(u) dx
\leq \frac{1}{2} \left(\frac{3t}{8\nu}+1 \right)^{1/2}
\int \theta^{1/2} |\nabla u| |u|^2 dx \ . $$
If $\gamma \geq 4$, we can use Young's inequality to conclude
$$\frac{1}{2} \left(\frac{3t}{8\nu}+1 \right)^{1/2}
\int \theta^{1/2} |\nabla u| |u|^2 dx \leq \frac{\nu}{2} \int \theta |\nabla u|^2 dx
+ \frac{1}{8\nu} \left(\frac{3t}{8\nu}+1 \right) \int |u|^4 dx \ . $$
In this case, $\| u \|_{L^4}^4 \leq C \| u \|_{L^\gamma}^4$ for some fixed $C$.
Using \eqref{CTU-gamma} and integrating in time, we obtain \eqref{CTU-upper}.\\
\\
Now assume that $\gamma \in [3,4]$.
We then use Young's inequality to get
$$ \frac{1}{2} \left(\frac{3t}{8\nu}+1 \right)^{1/2} \int \theta^{1/2} |\nabla u| |u|^2 dx \leq
\left( \frac{3t+8\nu}{32\nu} \right)^{1/2} \left( \int \theta
|\nabla u |^2 |u|^{\gamma-2} dx + \int |u|^{6-\gamma} dx \right) \nonumber \ . $$
We see that the first term is time-integrable by \eqref{CTU-gamma}. Since $\gamma \geq 3$,
we have that
$$ \int_{\TT^3} |u|^{6-\gamma} dx \leq C \| u \|_{L^\gamma}^{6-\gamma} \ , $$
for some fixed $C>0$. Integrating on $[0,T]$ and taking a supremum in time then yields
\begin{equation}
\sup_{t \in [0,T]} \left( \frac{1}{2} \| u (\cdot, t) \|_{L^2}^2 +
\| \theta(\cdot, t) \|_{L^2}^2 \right) \leq
\frac{1}{2} \| u_0 \|_{L^2}^2 + \| \theta_0 \|_{L^2}^2 +
C(T+1) \| u_0 \|_{L^\gamma}^{6-\gamma} \ , \label{CTU-upper-semi}
\end{equation}
for a different fixed constant $C>0$ that depends on $\nu$ and $\gamma$. Note that,
for $3 \leq \gamma \leq 4$, $6-\gamma < 3 \gamma - 3$.\\
\\
Lastly, we assume that $\gamma \in [9/4, 3)$ and employ the third term of \eqref{CTU-gamma}. Recalling
the Sobolev embedding of $H^2(\TT^3)$ into $L^6(\TT^3)$, there is some constant $C > 0$
such that
$$ \| u \|_{L^{3\gamma}}^\gamma = \left( \int |u|^{3\gamma} dx \right)^{1/3} \leq
C \int | \nabla |u|^{\gamma / 2}|^2 dx + C \| u \|_{L^\gamma}^\gamma \leq
C \frac{3t+8\nu}{8\nu} \int \theta | \nabla |u|^{\gamma / 2}|^2 dx +
C \| u \|_{L^\gamma}^\gamma \ . $$
For $ 3 > \gamma > 3/2$, we have that $\gamma < 6-\gamma < 3 \gamma$. Thus we may interpolate
$$ \| u \|_{L^{6-\gamma}} \leq \| u \|_{3 \gamma}^\alpha \| u \|_{L^\gamma}^{1-\alpha} \ , $$
where $\alpha = (9 -3\gamma)/(6-\gamma)$. Finally, since we assume $\gamma \geq 9/4$, we
have that $\alpha (6-\gamma) = 9-3\gamma \leq \gamma $. Therefore,
\begin{align}
\int |u|^{6-\gamma} dx &\leq C \| u \|_{L^\gamma}^{2\gamma - 3}
\left( \frac{3t+8\nu}{8\nu} \int \theta | \nabla |u|^{\gamma / 2}|^2 dx +
\| u \|_{L^\gamma}^\gamma \right)^{9/\gamma - 3} \nonumber \\
& \leq C \| u \|_{L^\gamma}^{2\gamma-3} \left( \frac{3t+8\nu}{8\nu}
\int \theta | \nabla |u|^{\gamma/2}|^2 dx
+ \| u \|_{L^\gamma}^\gamma +1 \right) \nonumber
\end{align}
Integrating on $[0,T]$ and taking a supremum in time then yields
$$ \sup_{t \in [0,T]} \left( \frac{1}{2} \| u (\cdot, t) \|_{L^2}^2 +
\| \theta(\cdot, t) \|_{L^2}^2 \right) \leq
\frac{1}{2} \| u_0 \|_{L^2}^2 + \| \theta_0 \|_{L^2}^2 + C(T^2+1)
(\| u_0 \|_{L^\gamma}^{3\gamma-3}+1) \ . $$
This, together with \eqref{CTU-upper-semi} yields \eqref{CTU-upper} .\\
\indent \hfill $\bullet$\\
\\
The proof of the weighted enstrophy bound for the thermal reduced Burgers equation will be very
similar to that for the thermal Navier-Stokes equations. Nevertheless, we believe it
is worth investigating since it is interesting to see how the lack of a divergence-free
condition on $u$ does not prevent the inequality from closing, despite being
an energy-type argument.\\
\\
\noindent{\sc Proof of Theorem 1.2:} Although
the statement of Theorem \ref{Thm-main} provides the thermal weight explicitly, we will
again perform the calculations with a generic weight $f(\theta)$. There will come
a point where $f$ will have to satisfy a specific ODE in order for the estimate to
proceed; at that point we will restrict to the weight mentioned in the theorem.
We do this to emphasize how (and exactly where) the methodology determines $f$.\\
\\
\noindent{}We start by integrating \eqref{CTU-F} against
$-\text{div}(f \nabla u)$ to get:
$$  \underbrace{\int f \frac{\partial_t}{2} |\nabla u|^2 dx}_{I_0} \underbrace{-
\int u_i \partial_i u_j \partial_k (f \partial_k u_j) dx - \frac{1}{\gamma} \int \text{div}(u)
u_j \partial_k (f \partial_k u_j) dx}_{I_A} \underbrace{+ \nu \int
\partial_i (\theta \partial_i u_j) \partial_k (f \partial_k u_j) dx}_{I_D} = 0 
$$
\\
The advection terms are straightforward:
\begin{align}
I_A &= \int f \partial_k u_i \partial_i u_j \partial_k u_j dx +
\int f u_i \partial_i \left( \frac{|\nabla u|^2}{2} \right) dx
+ \frac{1}{\gamma} \int f \text{div}(u) |\nabla u |^2 dx +
\frac{1}{\gamma} \int f \partial_k \text{div}(u) \partial_k u_j u_j dx \nonumber \\
&\leq - \frac{3\gamma - 2}{2\gamma} \int f |\nabla u|^3 dx -
\frac{1}{2} \int u_i \partial_i f |\nabla u|^2 dx
- \frac{1}{\gamma} \int f |\nabla^2 u| |\nabla u| |u| dx
\label{IA}
\end{align}
\\
Now we examine $I_D$, which will provide three positive terms (and one
cancellation) for our estimates.
\begin{align}
I_D &= \nu \int \partial_k(\theta \partial_i u_j) \partial_i (f \partial_k u_j)
dx = \nu \int (\partial_k \theta \partial_i u_j + \theta \partial_k \partial_i
u_j)(f' \partial_i \theta \partial_k u_j + f \partial_i \partial_k u_j) dx 
\nonumber \\
&= \nu \int \left( f' |\nabla \theta \cdot \nabla u|^2 + f \partial_k \theta
\partial_i u_j \partial_i \partial_k u_j + \theta f' \partial_i \theta
\partial_k u_j \partial_k \partial_i u_j + \theta f |\nabla^2 u|^2 \right) dx
\nonumber \\
&= \underbrace{\nu \int f' |\nabla \theta \cdot \nabla u|^2  dx}_{J_1} +
\underbrace{\nu \int (f + \theta f') \partial_k \theta \partial_i u_j \partial_i
\partial_k u_j dx}_{J_2} + \underbrace{\nu \int \theta f |\nabla^2 u|^2
dx}_{J_3} \nonumber
\end{align}
The last equality was obtained by re-indexing. Observe that $J_3$ is a positive
term and that $J_1$ has a sign depending on $f'$.
Continuing, we integrate $J_2$ by parts to move the $k$-derivative away
from $\partial_i u_j$ and onto $\partial_k \theta$:
\begin{align}
J_2 &= -\frac{\nu}{2} \int (f+\theta f') \partial_k \partial_k \theta |\nabla
u|^2 dx - \frac{\nu}{2} \int \partial_k (f+\theta f') \partial_k \theta |\nabla
u |^2 dx  \nonumber \\
&= \frac{\nu}{2 \kappa} \int \frac{f+\theta f'}{\theta} (-\kappa \theta \Delta \theta)
|\nabla u|^2 dx - \underbrace{ \frac{\nu}{2} \int (2f' + \theta f'') |\nabla \theta|^2
|\nabla u|^2 dx}_{K_0}  \nonumber \\
&= \frac{1}{2} \int \left( \frac{\nu}{\kappa} \frac{f+\theta f'}{\theta} \right)
\left(\underbrace{\frac{\nu}{2}|\nabla u|^2}_{K_1} \underbrace{+ 2 \kappa |
\nabla \theta|^2}_{K_2} \underbrace{- u \cdot \nabla \theta}_{K_3}
\underbrace{- \frac{\text{div}(u) \theta}{2}}_{K_4}
\underbrace{- \partial_t \theta}_{K_5} \right) |\nabla u|^2 dx - K_0  \nonumber
\end{align}
The last equality comes from direct substitution with \eqref{CTU-th}. We would
like to merge the integral involving the time derivative of $\theta$ (that is,
the one coming from the $K_5$ term) with $I_0$ to pull the time derivative
entirely outside the integral. In order to do this, we must clearly have
\begin{equation}
\frac{\nu}{\kappa} \frac{f(\theta)+\theta f'(\theta)}{\theta} = -f'(\theta)
\ \ \ \ \text{ or } \ \ \ \ f(\theta) = K \theta^{-\frac{1}{K}}
\ \ \ \ \text{ and } \ \ \ \ f'(\theta) = - \theta^{-\frac{K+1}{K}} \ . \nonumber
\end{equation}
Having established $f(\theta)$ explicitly, we now make some observations. As
mentioned before,
$$ I_0 + \frac{1}{2} \int \left( \frac{\nu}{\kappa} \frac{f+\theta f'}{\theta} \right) K_5
|\nabla u|^2 dx = \frac{\partial_t}{2} \int K \theta^{-\frac{1}{K}} |
\nabla u|^2 dx. $$
Moreover, the second term in \eqref{IA} from $I_A$ exactly cancels with the
integral arising from $K_3$.\\
\\
\noindent{}Also, $K_0$ and the integral arising from $K_2$ combine into
\begin{equation}
\nu \left( K - \frac{1}{2} - \frac{1}{2K} \right) \int \theta^{-\frac{K+1}{K}} |\nabla
\theta|^2 |\nabla u|^2 dx \label{KK}
\end{equation}
which is also positive.\\
\\
\noindent{}$J_1$ is again negative, but also controlled by
\eqref{KK} provided $K-1/(2K) > 3/2$.\\
\\
The term arising from $K_4$ is comparable to the first term in \eqref{IA}
from $I_A$. We therefore combine them.\\
\\
\noindent{}Collecting what terms remain, we obtain the inequality:
\begin{align}
\frac{\partial_t}{2} \int K \theta^{-\frac{1}{K}} & |\nabla u|^2 dx + \underbrace{\nu
K \int \theta^{\frac{K-1}{K}} |\nabla^2 u|^2 dx}_{U_1} +
\underbrace{\nu \frac{2K^2 - 3K - 1}{2K}
\int \theta^{-\frac{K+1}{K}} |\nabla \theta|^2 |\nabla u|^2 dx}_{U_2} \nonumber \\
&+ \underbrace{\frac{\nu}{4} \int \theta^{-\frac{K+1}{K}} |\nabla u|^4 dx}_{U_3}
\leq \underbrace{\left( \frac{3}{2}K+\frac{1}{4} \right)
\int \theta^{-\frac{1}{K}} |\nabla u|^3 dx}_{L_1}
+ \underbrace{\frac{K}{2} \int \theta^{-\frac{1}{K}}
|u| |\nabla u| |\nabla^2 u| dx}_{L_2}
\label{final-ineq}
\end{align}
\\
By H\"{o}lder's inequality,
\begin{equation}
L_1 \leq 2K \left(\int
|\nabla u|^6 dx \right)^{\frac{1}{9}} \left( \int |\nabla u|^4
\theta^{-\frac{K+1}{K}} dx \right)^{\frac{7}{12}} \left( \int \theta^2 dx
\right)^{\frac{7K-5}{24K}} | \TT^3 |^{\frac{7K-3}{24K}} \ . \nonumber
\nonumber
\end{equation}
Define
$$ \tau(t,\nu) := (3t/(8\nu) + 1) \ . $$
We then use Young's inequality with exponents $3$, $\frac{12}{7}$, and $12$ to obtain
\begin{equation}
L_1 \leq
\frac{\nu K}{4 \tau^{\frac{K-1}{K}}} \left( \int | \nabla u|^6 dx \right)^{\frac{1}{3}}
+ \frac{\nu}{12} \int \theta^{-\frac{K+1}{K}} | \nabla u |^4 dx
+ C \frac{K^8 \tau^{\frac{4(K-1)}{K}}}{\nu^{11}}
\| \theta \|_{L^2}^{\frac{7K-5}{K}} | \TT^3 |^{\frac{7K-3}{2K}}
\nonumber
\end{equation}
for some constant $C$ independent of $K$, $\nu$, and $\tau$,
and every $t \in [0,T]$. By \eqref{CTU-upper}, we have
$$ \int \theta^2 dx \leq \| \theta_0\|_{L^2}^2 + \frac{1}{2} \| u_0 \|_{L^2}^2 +
C(T^2 + 1)\left( \| u_0 \|_{L^\gamma}^{3\gamma-3} +1 \right) \ ,$$
and, by Sobolev embedding in three dimensions and the minimum principle \eqref{CTU-lower}
for $\theta$,
$$ \left( \int |\nabla u|^6 dx \right)^{\frac{1}{3}} \leq \int |\nabla^2 u|^2 dx
\leq \tau^{\frac{K-1}{K}} \int \theta^{\frac{K-1}{K}} |\nabla^2 u|^2 dx \ , $$
we see that $L_1$ of \eqref{final-ineq} is dominated by the
positive terms ($\frac{1}{4} U_1$ and $\frac{1}{3} U_3$) and the conserved
quantities of Proposition (\ref{Prop-initial}).\\
\\
\noindent{}Moving on to the last term, H\"{o}lder's inequality shows that
\begin{equation}
L_2 \leq \frac{K}{2}
\left( \int \theta^{\frac{K-1}{K}} | \nabla^2 u|^2 dx \right)^{\frac{1}{2}}
\left( \int \theta^{-\frac{K+1}{K}} | \nabla u |^4 dx \right)^{\frac{1}{4}}
\left( \int |u|^4 dx \right)^{\frac{1}{4}}
\| \theta^{-\frac{K+1}{4K}} \|_{L^\infty} \ .
\nonumber
\end{equation}
Then by Young's inequality and \eqref{CTU-lower},
\begin{equation}
L_2 \leq \frac{\nu K}{4} \int \theta^{\frac{K-1}{K}} | \nabla^2 u|^2 dx
+ \frac{\nu}{12} \int \theta^{-\frac{K+1}{K}} | \nabla u |^4 dx
+ C \frac{K}{\nu^3} \tau^{\frac{K+1}{4K}} \int |u|^4 dx \ .
\label{L_2-temp}
\end{equation}
The first and second terms are bounded by $\frac{1}{4} U_1$ and $\frac{1}{3} U_3$.
For the last term, if
$\gamma \geq 4$ we have $\| u \|_{L^4}^4 \leq C \| u \|_{L^\gamma}^4$, and we use
\eqref{CTU-gamma}. If $9/4 \leq \gamma < 4$ (or $\gamma=2$), we use interpolation
between $L^\gamma$ and $L^\infty$:
$$ \| u \|_{L^4} \leq \| u \|_{L^\gamma}^{\gamma / 4} \| u \|_{L^\infty}^{1-\gamma / 4} \ . $$
Agmon's inequality on a compact three-dimensional domain says that
$$ \| u \|_{L^\infty} \leq C \| u \|_{L^2}^{1/4} \| u \|_{H^2}^{3/4} \ . $$
We have from before that
$$ \| u \|_{\dot{H}^2}^2 \leq \tau^{\frac{K-1}{K}} \int \theta^{\frac{K-1}{K}}
|\nabla^2 u|^2 dx \ , $$
hence
$$ \| u \|_{L^4}^4 \leq C \| u \|_{L^\gamma}^\gamma \| u \|_{L^2}^{1-\gamma / 4}
\left( \tau^{\frac{K-1}{K}} \int \theta^{\frac{K-1}{K}} |\nabla^2 u|^2 dx +
\| u \|_{L^2}^2 \right)^{\frac{12-3 \gamma}{8}} \ . $$
Since $\gamma$ is in all cases larger than $4/3$, we have that $12-3\gamma < 8$.
So by Young's inequality again,
\begin{equation}
C \frac{K}{\nu^3} \tau^{\frac{K+1}{4K}} \| u \|_{L^4}^4 \leq
\frac{\nu K}{4} \int \theta^{\frac{K-1}{K}} |\nabla^2 u|^2 dx +
C \frac{K}{\nu^3} \tau^{\frac{K+1}{4K}}
\| u \|_{L^\gamma}^{\frac{8\gamma}{3\gamma - 4}} \| u \|_{L^2}^{\frac{8-2\gamma}{3\gamma-4}}
\left( \tau^{\frac{5K-3}{4K}} \nu^{-4} \right)^{\frac{12-3\gamma}{3\gamma-4}} \ ,
\label{L_2-last}
\end{equation}
where the $C$ on the right is a different constant from the one on the left.
Combining this with \eqref{L_2-temp}, we see that $L_2$ of \eqref{final-ineq} is
bounded similarly to $L_1$.\\
\\
Inequality \eqref{final-ineq} then reduces to
$$ \frac{K}{2} \partial_t \int \theta(x,t)^{-\frac{1}{K}} | \nabla u(x,t)|^2 dx
+\frac{1}{4} U_1 + U_2 + \frac{1}{3} U_3 \leq
C(\nu,\gamma,K, \| u_0 \|_{L^\gamma}, \| \theta_0 \|_{L^2})
\left( t^{M(K,\gamma)-1} + 1 \right) \ , $$
for some exponent $M(K,\gamma)$ which can be computed explicitly. Integrating this ODE
immediately gives \eqref{main-general} and \eqref{main-conserved}.\\
\indent \hfill $\bullet$


{\footnotesize

\begin{thebibliography}{100}
\bibitem{Agmon}
\textsc{S. Agmon}, {\it Lectures on elliptic boundary value problems} (2010), Providence,
RI: AMS Chelsea Publishing.\\
\bibitem{B-V}
\textsc{E. Baer and A. Vasseur}, {\it A bound from Below on the temperature for the
Navier-Stokes-Fourier System}, SIAM J. Math. Anal. \textbf{45} (2013), pp. 2046-2063.\\
\bibitem{BFM}
\textsc{M. Buli\"{c}ek, E. Feireisl, and J. M\'{a}lek}, {\it A Navier-Stokes-Fourier
system for incompressible fluids with temperature dependent material coefficients},
Nonlinear Analysis: Real World Applications \textbf{10} (2007), pp. 992-1015.\\
\bibitem{Vasseur}
\textsc{L. A. Caffarelli and A. F. Vasseur}, {\it The De Giorgi method for regularity of
solutions of elliptic equations and its applications to fluid dynamics}, Discrete and
Continuous Dynamical Systems \textbf{S 3} (2010), pp. 409-427.\\
\bibitem{Camp}
\textsc{S. Campanato}, {\it Propiet\`{a} di h\"{o}lderianit\`{a} di alcune classi di funzioni},
Ann. Scuola Norm. Sup. Pisa \textbf{17} (1963), pp. 175-188.\\
\bibitem{CRG}
\textsc{R. R. Coifman, R. Rochberg, and G. Weiss}, {\it Factorization theorems for Hardy spaces
in several variables}, Ann. of Math. \textbf{103} (1976), pp. 611-635.\\
\bibitem{D-F}
\textsc{B. Ducomet and E. Feireisl}, {\it On the Dynamics of Gaseous Stars}, Arch.
Rational Mech. Anal. \textbf{174} (2004), pp. 221-266.\\
\bibitem{Feireisl}
\textsc{E. Feireisl}, {\it Dynamics of Viscous Compressible Fluids} (2003),
Oxford Lecture Series in Mathematics and Its Applications.\\
\bibitem{Feireisl-again}
\textsc{E. Feireisl}, {\it On the motion of a viscous, compressible, and heat conducting
fluid}, Indiana Univ. Math. Journal \textbf{53} (2004), pp. 1705-1738.\\
\bibitem{F-N}
\textsc{E. Feireisl and A. Novotn\'{y}}, {\it Singular limits in thermodynamics of viscous
fluids}, Birkh\"{a}user-Verlag, Basel (2009).\\
\bibitem{F-M}
\textsc{E. Feireisl and J. M\'{a}lek}, {\it On the Navier-Stokes equations with
temperature-dependent transport coefficients}, Differ. Equ. Nonlinear Mech. (2006)
14pp. (electronic) Art. ID 90616.\\
\bibitem{HMW}
\textsc{R. A. Hunt, B. Muckenhoupt, and R. L. Wheeden}, {\it Weighted norm inequalities for the
conjugate function and the Hilbert transform}, Trans. Amer. Math. Soc. \textbf{176} (1973),
pp. 227-251.\\
\bibitem{JN}
\textsc{R. Johnson and C. J. Neugebauer}, {\it Properties of BMO functions whose reciprocals
are also BMO}, Z. Anal. Anwend. \textbf{12} (1993), pp. 3-11.\\
\bibitem{Lady}
\textsc{O. A. Ladyzhenskaya}, {\it Uniqueness and smoothness of generalized solutions of
Navier-Stokes equations}, Zap. Nauchn. Sem. Leningrad. Otdel. Mat. Inst. Steklov.
\textbf{5} (1967), pp. 169-185; English transl., Sem. Math. Steklov Math. Inst.,
Leningrad \textbf{5} (1969), pp. 60-67.\\
\bibitem{Lautrup}
\textsc{B. Lautrup}, {\it Physics of Continuous Matter} (2011), CRC Press Taylor and
Francis Group.\\
\bibitem{M-N}
\textsc{A. Matsumura and T. Nishida}, {\it Initial-boundary value problem for the equations of
motion of compressible viscous and heat-conductive fluids}, Comm. Math. Phys.
\textbf{89} (1983), pp. 445-464.\\
\bibitem{M-V}
\textsc{A. Mellet and A. Vasseur}, {\it A bound from below for the temperature in
compressible Navier-Stokes equations}, Monatsh. Math. \textbf{157} (2009), pp. 143-161.\\
\bibitem{Pet}
\textsc{S. Petermichl}, {\it The sharp weighted bound for the Riesz transforms},
Proc. Amer. Math. Soc. \textbf{136} (2007), pp. 1237-1249.\\
\bibitem{S-V}
\textsc{L. Silvestre and V. Vicol}, {\it H\"{o}lder continuity for a drift-diffusion
equation with pressure}, Annales de l'Institut Henri Poincar\'{e} (C) Analyse Non
Lin\'{e}aire \textbf{29} (2012), pp. 637-652.\\
\bibitem{Ser}
\textsc{J. Serrin}, {\it The initial value problem for the Navier-Stokes equations},
Nonlinear Problems, Univ. Wisconsin Press, Madison, 1963, pp. 69-98.\\
\bibitem{Stein}
\textsc{E. Stein}, {\it Singular Integrals and Differentiability Properties of Functions},
Princeton University Press, 1970.\\
\bibitem{Tao}
\textsc{T. Tao}, {\it Finite time blowup for an averaged three-dimensional Navier-Stokes
equation}, J. Amer. Math. Soc. \textbf{29} (2016), pp. 601-674.\\
\bibitem{Unter}
\textsc{J. Unterberger}, {\it Global existence for strong solutions of viscous Burgers
equation. (1) The bounded case}, (2015) arXiv:1503.05145.
\end{thebibliography}
}
\end{document}